\documentclass[12pt, a4paper]{amsart}
\usepackage[utf8]{inputenc}
\usepackage{amsxtra,amssymb,amsthm,amsmath,amscd,mathrsfs, epsfig, eufrak}
\usepackage{amscd, amsmath, mathrsfs, amssymb, amsthm, amsxtra, bbding, epsfig, eucal, eufrak, graphicx, latexsym, mathrsfs, mathbbol, bbold}
\usepackage[all]{xy}

\usepackage{tikz}

\usepackage[normalem]{ulem}
\usepackage{soul}
\usepackage{color}

\setstcolor{red}

\makeatletter
\@namedef{subjclassname@2010}{%
\textup{2010} Mathematics Subject Classification}
\makeatother

\def \P {{\mathbb P}}

\def \R {{\mathbb R}}

\def \d {\,{\rm d}}

\def\le{\leqslant}
\def\ge{\geqslant}

\theoremstyle{plain}
\newtheorem{theorem}{Theorem}
\newtheorem{proposition}{Proposition}[section]
\newtheorem{lemma}[proposition]{Lemma}

\theoremstyle{remark}

\numberwithin{equation}{section}

\addtocounter{footnote}{1}

\begin{document}


\title[On some sums involving the integral part function]
{On some sums involving the integral part function}
\author[\tiny Kui Liu, Jie Wu \& Zhishan Yang]{Kui Liu, Jie Wu \& Zhishan Yang}

\address{%
Kui Liu
\\
School of Mathematics and Statistics
\\
Qingdao University
\\
308 Ningxia Road
\\
Qingdao
\\
Shandong 266071
\\
China}
\email{liukui@qdu.edu.cn}

\address{%
Jie Wu
\\
CNRS UMR 8050\\
Laboratoire d'Analyse et de Math\'ematiques Appliqu\'ees\\
Universit\'e Paris-Est Cr\'eteil\\
94010 Cr\'eteil cedex\\
France}
\email{jie.wu@u-pec.fr}

\address{%
Zhishan Yang
\\
School of Mathematics and Statistics
\\
Qingdao University
\\
308 Ningxia Road
\\
Qingdao
\\
Shandong 266071
\\
China}
\email{zsyang@qdu.edu.cn}

\date{\today}

\begin{abstract}
Denote by $\tau_k(n)$, $\omega(n)$ and $\mu_2(n)$ 
the number of representations of $n$ as product of $k$ natural numbers, 
the number of distinct prime factors of $n$ 
and the characteristic function of the square-free integers, respectively.
Let $[t]$ be the integral part of real number $t$.
For $f=\omega, 2^{\omega}, \mu_2, \tau_k$, we prove that
$$
\sum_{n\le x} f\Big(\Big[\frac{x}{n}\Big]\Big)
= x\sum_{d\ge 1} \frac{f(d)}{d(d+1)} + O_{\varepsilon}(x^{\theta_f+\varepsilon})
$$
for $x\to\infty$, where 
$\theta_{\omega}=\frac{53}{110}$,
$\theta_{2^{\omega}}=\frac{9}{19}$,
$\theta_{\mu_2}=\frac{2}{5}$,
$\theta_{\tau_k}=\frac{5k-1}{10k-1}$
and $\varepsilon>0$ is an arbitrarily small positive number.
These improve the corresponding results of Bordell\`es.
\end{abstract}


\maketitle

\section{Introduction}

Denote by $[t]$ the integral part of the real number $t$.
Recently Bordell\`es, Dai, Heyman, Pan and Shparlinski \cite{BDHPS2019} established an asymptotic formula of
$$
S_f(x) := \sum_{n\le x} f\Big(\Big[\frac{x}{n}\Big]\Big)
$$
under some simple assumptions of $f$.
Subsequently, Wu \cite{Wu2020} and Zhai \cite{Zhai2020} improved their results independently.
In particular, if $f(n)\ll_{\varepsilon} n^{\varepsilon}$ for any $\varepsilon>0$ and all $n\ge 1$,
then \cite[Theorem 1.2(i)]{Wu2020} or \cite [Theorem 1]{Zhai2020} yields the following asymptotic formula
\begin{equation}\label{WuZhai}
S_f(x) = x\sum_{d\ge 1}\frac{f(d)}{d(d+1)} + O_{\varepsilon}(x^{1/2+\varepsilon})
\end{equation}
as $x\to\infty$.
Ma and Wu \cite{MaWu2020} observed that if $f$ is factorizable in certain sense, 
then it is possible to break the $\frac{1}{2}$-barrier for the error term of \eqref{WuZhai}.
In particular, for the von Mangoldt function $\Lambda(n)$,
they proved, with the help of the Vaughan identity and the technique of one-dimensional exponential sum, that
\begin{equation}\label{MW:Lambda}
S_{\Lambda}(x)
= x\sum_{d\ge 1}\frac{\Lambda(d)}{d(d+1)} + O_{\varepsilon}(x^{35/71+\varepsilon})
\qquad
(x\to\infty).
\end{equation}
Using similar idea to the classical divisor function $\tau(n)$, Ma and Sun \cite{MaSun2020} showed that
\begin{equation}\label{MS:tau}
S_{\tau}(x)
= x\sum_{d\ge 1}\frac{\tau(d)}{d(d+1)} + O_{\varepsilon}(x^{11/23+\varepsilon})
\qquad
(x\to\infty).
\end{equation}
Subsequently,
with the help of a result of Baker on 2-dimensional exponential sums \cite[Theorem 6]{Baker2007},
Bordell\`es \cite{Bordelles2020} sharpened the exponents in \eqref{MW:Lambda}-\eqref{MS:tau}
and also studied some new examples.
Denote by $\tau_k(n)$, $\omega(n)$ and $\mu_2(n)$ 
the number of representations of $n$ as product of $k$ natural numbers, 
the number of distinct prime factors of $n$ 
and the characteristic function of the square-free integers, respectively.
His results can be stated as follows:
\begin{equation}\label{Result:Bordelles}
\sum_{n\le x} f\Big(\Big[\frac{x}{n}\Big]\Big)
= x\sum_{d\ge 1} \frac{f(d)}{d(d+1)} + O_{\varepsilon}(x^{\vartheta_f+\varepsilon}),
\end{equation}
with
$$\displaylines{
\vbox{\tabskip = 0mm\offinterlineskip
\halign{
\vrule # & &\hfil$ $ $#$  \hfil & \vrule #\cr
\noalign{\hrule}
height 1,5mm
&&&&&&&&&&&&&&&&
\cr
& \;\; f \;\; & 
& \;\; \Lambda \;\; & 
& \;\; \omega \;\; & 
& \;\; 2^{\omega} \;\; & 
& \;\; \mu_2 \;\; & 
& \;\; \tau \;\; & 
& \;\; \tau_3 \;\; & 
& \;\; \tau_k \; {\scriptstyle{(k\,\ge\, 4)}} \;\; & 
\cr
height 1,5mm
&&&&&&&&&&&&&&&&
\cr
\noalign{\hrule}
height 1,5mm
&&&&&&&&&&&&&&&&
\cr
& \vartheta_f &
& \tfrac{97}{203} & 
& \tfrac{455}{914} &
& \tfrac{97}{202} &
& \tfrac{1919}{4268} &
& \tfrac{19}{40} &
& \tfrac{283}{574} &
& \tfrac{4k^3-k-2}{8k^3-2k-2} &
\cr
height 1,5mm
&&&&&&&&&&&&&&&&
\cr
\noalign {\hrule}}}
\cr}$$
Very recently, Liu, Wu and Yang \cite{LiuWuYang2021} succeeded to improve Bordell\`es' exponent 
$\frac{97}{203}$ to $\frac{9}{19}$ for $\Lambda$.
Their principal tool is a new estimate on 3-dimensional exponential sums.

The aim of this paper is to propose better exponents than all others cases of the above table.

\begin{theorem}\label{thm}
Under the previous notation, 
the asymptotic formula \eqref{Result:Bordelles} holds with
$$\displaylines{
\vbox{\tabskip = 0mm\offinterlineskip
\halign{
\vrule # & &\hfil$ $ $#$  \hfil & \vrule #\cr
\noalign{\hrule}
height 1,5mm
&&&&&&&&&&
\cr
& f & 
& \omega & 
& 2^{\omega} & 
& \mu_2 & 
& \tau_k \; {\scriptstyle{(k\,\ge\, 2)}} & 
\cr
height 1,5mm
&&&&&&&&&&
\cr
\noalign{\hrule}
height 1,5mm
&&&&&&&&&&
\cr
& \;\; \vartheta_f \;\; &
& \;\; \frac{53}{110} \;\; & 
& \;\; \frac{9}{19} \;\; & 
& \;\; \frac{2}{5} \;\; & 
& \;\; \frac{5k-1}{10k-1} \;\; &
\cr
height 1,5mm
&&&&&&&&&&
\cr
\noalign {\hrule}}}
\cr}$$
\end{theorem}

The key point of proving Theorem \ref{thm} is to establish good bounds for
$$
\mathfrak{S}_{\delta}^{f}(x, D)
:= \sum_{D<d\le 2D} f(d) \psi\Big(\frac{x}{d+\delta}\Big).
$$
where $\psi(t) := \{t\}-\frac{1}{2}$ and $\{t\}$ means the fractional part of real number $t$.
Let $\P$ be the set of all primes and define
\begin{equation}\label{def:functions}
\mathbb{1}_{\P}(n)
:= \begin{cases}
1 & \text{if $n\in \P$},
\\
0 & \text{otherwise},
\end{cases}
\qquad
\mathbb{1}(n)\equiv 1,
\qquad
\tilde{\mu}(n)
:= \begin{cases}
\mu(d) & \text{if $n=d^2$},
\\
0 & \text{otherwise}.
\end{cases}
\end{equation}
Then we have the following relations:
$$
\omega = \mathbb{1}_{\P}*\mathbb{1},
\qquad
\tau_k=\tau_{k-1}*\mathbb{1},
\qquad
2^{\omega} = \mu_2*\tau,
\qquad
\mu_2=\tilde{\mu}*\mathbb{1}.
$$
Thus when $f=\omega, \tau_k, 2^{\omega}, \mu_2$, we can use these relations
to decompose $\mathfrak{S}_{\delta}^{f}(x, D)$ into bilinear forms
\begin{align}
\mathfrak{S}_{\delta}^{\omega}(x, D)
& := \sum_{D<d\ell\le 2D} \mathbb{1}_{\P}(d) \psi\Big(\frac{x}{d\ell+\delta}\Big),
\label{def:Somega}
\\
\mathfrak{S}_{\delta}^{\tau_k}(x, D)
& := \sum_{D<d\ell\le 2D} \tau_{k-1}(d) \psi\Big(\frac{x}{d\ell+\delta}\Big),
\label{def:Stauk}
\\
\mathfrak{S}_{\delta}^{2^{\omega}}(x, D)
& := \sum_{D<d^2\ell\le 2D} \mu(d)\tau(\ell) \psi\Big(\frac{x}{d^2\ell+\delta}\Big),
\label{def:S2omega}
\\
\mathfrak{S}_{\delta}^{\mu_2}(x, D)
& := \sum_{D<d^2\ell\le 2D} \mu(d) \psi\Big(\frac{x}{d^2\ell+\delta}\Big).
\label{def:Smu2}
\end{align}
In the third section, we shall use the Fourier analyse and the technique of multiple exponential sums
to establish our bounds for these bilinear forms.

\vskip 8mm

\section{Preliminary lemmas}

In this section, we shall cite three lemmas, which will be needed in the next section.
The first one is \cite[Proposition 3.1]{LiuWuYang2021}.

\begin{lemma}\label{lem:LiuWuYang}
Let $\alpha>0$, $\beta>0$, $\gamma>0$ and $\delta\in \R$ be some constants. 
For $X>0$, $H\ge 1$, $M\ge 1$ and $N\ge 1$, define
\begin{equation}\label{def:SdeltaHMN}
S_{\delta}
= S_{\delta}(H, M, N)
:= \sum_{h\sim H} \sum_{m\sim M} \sum_{n\sim N} a_{h, m} b_n
{\rm e}\bigg(X\frac{M^{\beta}N^{\gamma}}{H^{\alpha}} \frac{h^{\alpha}}{m^{\beta}n^{\gamma}+\delta}\bigg),
\end{equation}
where ${\rm e}(t):={\rm e}^{2\pi{\rm i}t}$, 
the $a_{h, m}$ and $b_n$ are complex numbers such that $|a_{h, m}|\le 1$ and $|b_n|\le 1$,
and $m\sim M$ means that $M<m\le 2M$.
For any $\varepsilon>0$ we have
$$
S_{\delta}
\ll \big((X^{\kappa}H^{2+\kappa}M^{2+\kappa}N^{1+\kappa+\lambda})^{1/(2+2\kappa)}
+ HM N^{1/2}
+ (HM)^{1/2} N
+ X^{-1/2}HMN\big)X^{\varepsilon} 
$$
uniformly for $M\ge 1$, $N\ge 1$, $H\le M^{\beta-1}N^{\gamma}$ and $0\le \delta\le 1/\varepsilon$, 
where $(\kappa, \lambda)$ is an exponent pair and 
the implied constant depends on $(\alpha, \beta, \gamma, \varepsilon)$ only.
\end{lemma}

\goodbreak

The second one is the Vaughan identity \cite[formula (3)]{Vaughan1980}.

\begin{lemma}\label{lem:VaughanIdentity}
There are six real arithmetic functions $\alpha_k(n)$ verifying
$\alpha_k(n)\ll \tau(n)\log(2n)$
for
$(n\ge 1, \, 1\le k\le 6)
$
such that for $D\ge 1$ and any arithmetical function $g$, we have
\begin{equation}\label{identity:Vaughan}
\sum_{D<d\le 2D} \Lambda(d) g(d)
= S_1 + S_2 + S_3 + S_4,
\end{equation}
where $\tau(n)$ is the classic divisor function and
\begin{align*}
S_1
& := \sum_{m\le D^{1/3}} \alpha_1(m) \sum_{D<\ell mn\le 2D} g(mn),
\\
S_2
& := \sum_{m\le D^{1/3}} \alpha_2(m) \sum_{D<\ell mn\le 2D} g(mn)\log n,
\\\noalign{\vskip 1mm}
S_3
& := \mathop{{\sum}\,\,{\sum}}_{\substack{D^{1/3}<m, n\le D^{2/3}\\ D<\ell mn\le 2D}} \alpha_3(m) \alpha_4(n) g(mn),
\\
S_4
& := \mathop{{\sum}\,\,{\sum}}_{\substack{D^{1/3}<m, n\le D^{2/3}\\ D<\ell mn\le 2D}} \alpha_5(m) \alpha_6(n) g(mn).
\end{align*}
The same result also holds for the M\"obius function $\mu$.
\end{lemma}

\vskip 1mm

The third one is due to Vaaler (see \cite[Theorem A.6]{GrahamKolesnik1991}).

\begin{lemma}\label{lem:Vaaler}
Let $\psi(t) := \{t\}-\frac{1}{2}$, where $\{t\}$ means the fractional part of real number $t$.
For $x\ge 1$ and $H\ge 1$, we have
$$
\psi(x) = - \sum_{1\le |h|\le H} \Phi\Big(\frac{h}{H+1}\Big) \frac{{\rm e}(hx)}{2\pi{\rm i} h} + R_H(x),
$$
where $\Phi(t):=\pi t (1-|t|)\cot(\pi t) + |t|$ and the second term $R_H(x)$ satisfies
\begin{equation}\label{eq:Vaaler}
|R_H(x)|\le \frac{1}{2H+2} \sum_{0\le |h|\le H} \Big(1-\frac{|h|}{H+1}\Big) {\rm e}(hx).
\end{equation}
\end{lemma}

\vskip 5mm

\section{bilinear forms}

The aim of this section is to establish some non-trivial bounds for the bilinear forms given by
\eqref{def:Somega}, \eqref{def:Smu2}, \eqref{def:Stauk} and \eqref{def:S2omega}.
For the convenience, we consider 
\begin{equation}\label{def:SLambda*1}
\mathfrak{S}_{\delta}^{\Lambda*\mathbb{1}}(x, D)
:= \sum_{D<d\ell\le 2D} \Lambda(d) \psi\Big(\frac{x}{d\ell+\delta}\Big),
\end{equation}
instead of $\mathfrak{S}_{\delta}^{\omega}(x, D)$.
The following proposition gives the required estimate for $\mathfrak{S}_{\delta}^{\Lambda*\mathbb{1}}(x, D)$.

\begin{proposition}\label{prop_1}
For any $\varepsilon>0$, we have
\begin{equation}\label{eq:prop_1}
\mathfrak{S}_{\delta}^{\Lambda*\mathbb{1}}(x, D)
\ll_{\varepsilon} \big((x^{19} D^{87})^{1/133}
+ (x^{171} D^{624})^{1/1026}
\big)x^{\varepsilon}
\end{equation}
uniformly for $x\ge 3$, $D\le x^{57/103}$ and $0\le \delta\le \varepsilon^{-1}$,
where the implied constant depends on $\varepsilon$ only.
\end{proposition}

\begin{proof}
Applying the Vaughan identity \eqref{identity:Vaughan} with $g(d) = \psi\big(\frac{x}{d\ell+\delta}\big)$, it follows that
\begin{equation}\label{Decomposition}
\mathfrak{S}_{\delta}^{\Lambda*\mathbb{1}}(x, D)
= \mathfrak{S}_{\delta, 1}^{\Lambda*\mathbb{1}} 
+ \mathfrak{S}_{\delta, 2}^{\Lambda*\mathbb{1}} 
+ \mathfrak{S}_{\delta, 3}^{\Lambda*\mathbb{1}} 
+ \mathfrak{S}_{\delta, 4}^{\Lambda*\mathbb{1}},
\end{equation}
where
\begin{align*}
\mathfrak{S}_{\delta, 1}^{\Lambda*\mathbb{1}}
& := \sum_{\ell\le D}
\sum_{m\le (D/\ell)^{1/3}} \alpha_1(m) \sum_{D<\ell mn\le 2D} \psi\Big(\frac{x}{\ell mn+\delta}\Big),
\\\noalign{\vskip 1mm}
\mathfrak{S}_{\delta, 2}^{\Lambda*\mathbb{1}}
& := \sum_{\ell\le D}
\sum_{m\le (D/\ell)^{1/3}} \alpha_2(m) \sum_{D<\ell mn\le 2D} \psi\Big(\frac{x}{\ell mn+\delta}\Big)\log n,
\\\noalign{\vskip 1mm}
\mathfrak{S}_{\delta, 3}^{\Lambda*\mathbb{1}}
& := \sum_{\ell\le D}
\mathop{{\sum}\,\,{\sum}}_{\substack{(D/\ell)^{1/3}<m, n\le (D/\ell)^{2/3}\\ D<\ell mn\le 2D}}
\alpha_3(m) \alpha_4(n) \psi\Big(\frac{x}{\ell mn+\delta}\Big),
\\\noalign{\vskip -1mm}
\mathfrak{S}_{\delta, 4}^{\Lambda*\mathbb{1}}
& := \sum_{\ell\le D}
\mathop{{\sum}\,\,{\sum}}_{\substack{(D/\ell)^{1/3}<m, n\le (D/\ell)^{2/3}\\ D<\ell mn\le 2D}}
\alpha_5(m) \alpha_6(n) \psi\Big(\frac{x}{\ell mn+\delta}\Big).
\end{align*}
Using Lemma \ref{lem:Vaaler} and splitting the intervals of summation into the dyadic intervals, we can write
\begin{equation}\label{proof:prop_2}
\begin{aligned}
\mathfrak{S}_{\delta, j}^{\Lambda*\mathbb{1}}
& = - \frac{1}{2\pi\text{i}} \sum_{H'} \sum_{L} \sum_{M} \sum_{N} 
\big(\mathfrak{S}_{\delta, j, \flat}^{\Lambda*\mathbb{1}}(H', L, M, N) 
+ \overline{\mathfrak{S}_{\delta, j, \flat}^{\Lambda*\mathbb{1}}(H', L, M, N)}\big)
\\
& \hskip 31mm
+ \sum_{L} \sum_{M} \sum_{N} \mathfrak{S}_{\delta, j, \dagger}^{\Lambda*\mathbb{1}}(L, M, N),
\end{aligned}
\end{equation}
where $1\le H'\le H\le D$, $a_h:=\frac{H'}{h}\Phi\big(\frac{h}{H+1}\big)\ll 1$ and 
\begin{align*}
\mathfrak{S}_{\delta, 1, \flat}^{\Lambda*\mathbb{1}}
= \mathfrak{S}_{\delta, 1, \flat}^{\Lambda*\mathbb{1}}(H', L, M, N)
& := \frac{1}{H'} \sum_{h\sim H'} \mathop{\sum_{\ell\sim L} \sum_{m\sim M} \sum_{n\sim N}}_{D<\ell mn\le 2D} 
a_h\alpha_1(m) \text{e}\Big(\frac{hx}{\ell mn+\delta}\Big),
\\
\mathfrak{S}_{\delta, 1, \dagger}^{\Lambda*\mathbb{1}}
= \mathfrak{S}_{\delta, 1, \dagger}^{\Lambda*\mathbb{1}}(L, M, N)
& := \mathop{\sum_{\ell\sim L} \sum_{m\sim M} \sum_{n\sim N}}_{D<\ell mn\le 2D} \alpha_1(m) 
R_H\Big(\frac{x}{\ell mn+\delta}\Big),
\end{align*}
and $\mathfrak{S}_{\delta, j, \flat}^{\Lambda*\mathbb{1}}(H', L, M, N)$,
$\mathfrak{S}_{\delta, j, \dagger}^{\Lambda*\mathbb{1}}(L, M, N)$ can be defined similarly for $j=2, 3, 4$.
We have
\begin{equation}\label{Cond:LMN-Type1}
LMN\asymp D,
\qquad
M\le (D/L)^{1/3}
\quad\text{for}\quad
j=1, 2
\end{equation}
and
\begin{equation}\label{Cond:LMN-Type2}
LMN\asymp D,
\quad
(D/L)^{1/3}\le M\le (D/L)^{1/2}\le N\le (D/L)^{2/3} 
\;\;\text{for}\;\;
j=3, 4.
\end{equation}
In the last case, we have considered the symmetry of $m$ and $n$.

\vskip 1,5mm
\goodbreak

A. \textit{Bounds of $\mathfrak{S}_{\delta, j}^{\Lambda*\mathbb{1}}$ for $j=1, 2$.}
\par
\vskip 0,5mm

Since we can remove the factor $\log n$ by a simple partial integration, 
we only treat $\mathfrak{S}_{\delta, 1}^{\Lambda*\mathbb{1}}$.
If $L\le D^{1/5}$, applying the exponent pair $(\frac{1}{6}, \frac{4}{6})$ to the sum over $n$, we have
\begin{align*}
\mathfrak{S}_{\delta, 1, \flat}^{\Lambda*\mathbb{1}}
& \ll x^{\varepsilon} \frac{1}{H'}
\sum_{h\sim H'} \sum_{\ell\sim L} \sum_{m\sim M}
\bigg(\Big(\frac{xH'}{LMN^2}\Big)^{1/6} N^{4/6} 
+ \Big(\frac{xH'}{LMN^2}\Big)^{-1}\bigg)
\\
& \ll x^{\varepsilon} 
\big((xHL^5M^5N^2)^{1/6} 
+ x^{-1}(LMN)^2\big).
\end{align*}
In view of $L\le D^{1/5}$ and \eqref{Cond:LMN-Type1}, we derive that 
\begin{equation}\label{S1flat:1}
\mathfrak{S}_{\delta, 1, \flat}^{\Lambda*\mathbb{1}}
\ll \big((x^5D^{17}H^5)^{1/30} + x^{-1}D^2\big)x^{\varepsilon}
\end{equation}
for $H'\le H\le D$ and $(L, M, N)$ verifying \eqref{Cond:LMN-Type1}.

Next we suppose that $D^{1/5}\le L\le D^{9/19}$.
Firstly we remove the extra multiplicative condition $D<\ell mn\le 2D$ at the cost of a factor $(\log x)^2$
and then apply Lemma \ref{lem:LiuWuYang} with $\alpha=\beta=\gamma=1$,
$(X, H, M, N) = (xH'/LMN, H', LM, N)$ and $(\kappa, \lambda)=(\frac{1}{2}, \frac{1}{2})$ to get
\begin{equation}\label{S1flat:2}
\begin{aligned}
\mathfrak{S}_{\delta, 1, \flat}^{\Lambda*\mathbb{1}}
& \ll \big((x L^4M^4 N^3)^{1/6}
+ LM N^{1/2}
+ (LM)^{1/2} N
+ x^{-1/2}(LMN)^{3/2}
\big) x^{\varepsilon}
\\
& \ll \big((x^3 D^{10} L^2)^{1/18}
+ (D^2L)^{1/3}
+ DL^{-1/2}
+ (x^{-1}D^3)^{1/2}
\big) x^{\varepsilon}
\\
& \ll \big((x^{171} D^{624})^{1/1026}
+ D^{9/10}
+ (x^{-1}D^3)^{1/2}
\big) x^{\varepsilon}
\end{aligned}
\end{equation}
for $H'\le H\le D^{20/57} \; (\le (D/L)^{2/3}\le N$)
and $(L, M, N)$ verifying \eqref{Cond:LMN-Type1} and $D^{1/5}\le L\le D^{9/19}$.

When $D^{9/19}\le L\le D$, we apply the exponent pair $(\frac{1}{6}, \frac{4}{6})$ to the sum over $\ell$ to get
\begin{align*}
\mathfrak{S}_{\delta, 1, \flat}^{\Lambda*\mathbb{1}}
& \ll x^{\varepsilon} \frac{1}{H'} \sum_{h\sim H'} \sum_{m\sim M} \sum_{n\sim N}
\bigg(\Big(\frac{xH'}{L^2MN}\Big)^{1/6} L^{4/6} + \Big(\frac{xH'}{L^2MN}\Big)^{-1}\bigg)
\\
& \ll x^{\varepsilon}
\big((xHL^2M^5N^5)^{1/6} + x^{-1} H'^{-1}(LMN)^2\big).
\end{align*}
From this we deduce that
\begin{equation}\label{S1flat:3}
\mathfrak{S}_{\delta, 1, \flat}^{\Lambda*\mathbb{1}}
\ll \big((x^{19} D^{68} H^{19})^{1/114} + x^{-1} D^2\big)x^{\varepsilon}
\end{equation}
for $H'\le H\le D$ and $(L, M, N)$ verifying \eqref{Cond:LMN-Type1} and $D^{9/19}\le L\le D$.

Combining \eqref{S1flat:1}, \eqref{S1flat:2} and \eqref{S1flat:3}, we obtain
\begin{equation}\label{S1flat:4}
\begin{aligned}
\mathfrak{S}_{\delta, 1, \flat}^{\Lambda*\mathbb{1}}
& \ll \big((x^{171} D^{624})^{1/1026}
+ D^{9/10}
+ (x^{-1}D^3)^{1/2}
\\
& \hskip 3,5mm
+ (x^{19} D^{68} H^{19})^{1/114}
+ (x^5D^{17}H^5)^{1/30}
\big)x^{\varepsilon}
\end{aligned}
\end{equation}
for $H'\le H\le D^{20/57}$ and $(L, M, N)$ verifying \eqref{Cond:LMN-Type1},
where we have used the fact that $x^{-1}D^2\le (x^{-1}D^3)^{1/2}$.
The same bound holds for $\mathfrak{S}_{\delta, 2}^{\Lambda, \dagger}$.
Therefore for $j=1, 2$, we have
\begin{equation}\label{S1S2}
\begin{aligned}
\mathfrak{S}_{\delta, j}^{\Lambda*\mathbb{1}}
& \ll \big((x^{171} D^{624})^{1/1026}
+ D^{9/10}
+ (x^{-1}D^3)^{1/2}
\\
& \hskip 3,5mm
+ (x^{19} D^{68} H^{19})^{1/114}
+ (x^5D^{17}H^5)^{1/30}
\big)x^{\varepsilon}
\end{aligned}
\end{equation}
for $H\le D^{20/57}$ and $(L, M, N)$ verifying \eqref{Cond:LMN-Type1}.

\vskip 1,5mm
\goodbreak

B. \textit{Bounds of $\mathfrak{S}_{\delta, j}^{\Lambda*\mathbb{1}}$ for $j=3, 4$.}
\par
\vskip 0,5mm

In this case, we should estimate
\begin{align*}
\mathfrak{S}_{\delta, 3, \flat}^{\Lambda*\mathbb{1}}
= \mathfrak{S}_{\delta, 3, \flat}^{\Lambda*\mathbb{1}}(H', L, M, N)
& := \frac{1}{H'} \sum_{h\sim H'} \mathop{\sum_{\ell\sim L} \sum_{m\sim M} \sum_{n\sim N}}_{D<\ell mn\le 2D} 
a_h\alpha_3(m) \alpha_4(n) \text{e}\Big(\frac{hx}{\ell mn+\delta}\Big),
\\
\mathfrak{S}_{\delta, 3, \dagger}^{\Lambda*\mathbb{1}}
= \mathfrak{S}_{\delta, 3, \dagger}^{\Lambda*\mathbb{1}}(L, M, N)
& := \mathop{\sum_{\ell\sim L} \sum_{m\sim M} \sum_{n\sim N}}_{D<\ell mn\le 2D} \alpha_3(m) \alpha_4(n) 
R_H\Big(\frac{x}{\ell mn+\delta}\Big),
\end{align*}
for $(L, M, N)$ verifying \eqref{Cond:LMN-Type2}.

Firstly we bound $\mathfrak{S}_{\delta, 3, \flat}^{\Lambda*\mathbb{1}}(H', L, M, N)$.
We remove the extra multiplicative condition $D/\ell<mn\le 2D/\ell$ at the cost of a factor $\log x$.
By Lemma \ref{lem:LiuWuYang} with $\alpha=\beta=\gamma=1$, $(X, H, M, N)=(xH'/LMN, H', M, LN)$
and $(\kappa, \lambda)=(\frac{1}{2}, \frac{1}{2})$, we can derive 
\begin{equation}\label{S3flat:1}
\begin{aligned}
\mathfrak{S}_{\delta, 3, \flat}^{\Lambda*\mathbb{1}}
& \ll_{\varepsilon}\big((x L^3 M^4 N^3)^{1/6}
+ M (LN)^{1/2}
+ M^{1/2} LN
+ (x^{-1}DH')^{1/2}\big) x^{\varepsilon}
\\
& \ll_{\varepsilon} \big((x^2 D^7 L^{-1})^{1/12}
+ (D^5L)^{1/6}
+ (x^{-1}DH')^{1/2}\big) x^{\varepsilon},
\end{aligned}
\end{equation}
for $H'\le H\le D^{1/2} \; (\le (LD)^{1/2}\le LN)$ and $(L, M, N)$ verifying \eqref{Cond:LMN-Type2},
where we have removed the term $(D^3 L^{-1})^{1/4}\;(\le (D^5L)^{1/6})$.
On the other hand, we can apply the exponent pair $(\frac{1}{6}, \frac{4}{6})$ to the sum over $\ell$ 
(similar to \eqref{S1flat:1}) to get
\begin{equation}\label{S3flat:2}
\mathfrak{S}_{\delta, 3, \flat}^{\Lambda*\mathbb{1}}
\ll \big((xD^5HL^{-3})^{1/6} + x^{-1} D^2\big)x^{\varepsilon}
\end{equation}
for $H'\le H\le D$ and $(L, M, N)$ verifying \eqref{Cond:LMN-Type2}.

From \eqref{S3flat:1} and \eqref{S3flat:2}, we deduce that
\begin{equation}\label{S3flat:3}
\begin{aligned}
\mathfrak{S}_{\delta, 3, \flat}^{\Lambda*\mathbb{1}}
& \ll_{\varepsilon} \big((x^2 D^7 L^{-1})^{1/12}
+ (xD^{20}H)^{1/24}
+ (x^{-1}DH)^{1/2}
+ x^{-1} D^2\big) x^{\varepsilon},
\end{aligned}
\end{equation}
for $H\le D^{1/2}$ and $(L, M, N)$ verifying \eqref{Cond:LMN-Type2}, 
where we have used the following estimate
$$
\min\{(D^5L)^{1/6}, \, (xD^5HL^{-3})^{1/6}\}
\le ((D^5L)^{1/6})^{\frac{18}{24}}
((xD^5HL^{-3})^{1/6})^{\frac{6}{24}} 
= (xD^{20}H)^{1/24}.
$$

Secondly we bound $\mathfrak{S}_{\delta, 3, \dagger}^{\Lambda*\mathbb{1}}$.
Using \eqref{eq:Vaaler} of Lemma \ref{lem:Vaaler}, we have
\begin{align*}
\mathfrak{S}_{\delta, 3, \dagger}^{\Lambda*\mathbb{1}}
& \ll x^{\varepsilon} \sum_{\ell\sim L} \sum_{m\sim M} \sum_{n\sim N}
\Big|R_H\Big(\frac{x}{\ell mn+\delta}\Big)\Big|
\\\noalign{\vskip -0,5mm}
& \ll \frac{x^{\varepsilon}}{H} \sum_{\ell\sim L} \sum_{m\sim M} \sum_{n\sim N}
\sum_{0\le |h|\le H} \Big(1-\frac{|h|}{H+1}\Big) {\rm e}\Big(\frac{hx}{\ell mn+\delta}\Big)
\\\noalign{\vskip -0,5mm}
& \ll x^{\varepsilon}\big(DH^{-1} 
+ \max_{1\le H'\le H} \big|\widetilde{\mathfrak{S}}_{\delta, 3, \dagger}^{\Lambda*\mathbb{1}}(H', L, M, N)\big|\big),
\end{align*}
where
$$
\widetilde{\mathfrak{S}}_{\delta, 3, \dagger}^{\Lambda*\mathbb{1}}(H', L, M, N)
:= \frac{1}{H} \sum_{\ell\sim L} \sum_{m\sim M} \sum_{n\sim N} \sum_{h\sim H'} 
\Big(1-\frac{|h|}{H+1}\Big) {\rm e}\Big(\frac{hx}{\ell mn+\delta}\Big).
$$
Clearly we can bound $\widetilde{\mathfrak{S}}_{\delta, 3, \dagger}^{\Lambda*\mathbb{1}}(H', L, M, N)$
in the same way as $\mathfrak{S}_{\delta, 3, \flat}^{\Lambda*\mathbb{1}}(H', L, M, N)$
and get
\begin{equation}\label{S3dagger:1}
\begin{aligned}
\mathfrak{S}_{\delta, 3, \dagger}^{\Lambda*\mathbb{1}}
& \ll_{\varepsilon} 
\big(DH^{-1}
+ (x^2 D^7 L^{-1})^{1/12}
+ (xD^{20}H)^{1/24}
+ (x^{-1}DH)^{1/2}
+ x^{-1} D^2
\big) x^{\varepsilon}
\end{aligned}
\end{equation}
for $H\le D^{1/2}$ and $(L, M, N)$ verifying \eqref{Cond:LMN-Type2}.
Thus \eqref{S3flat:3} and \eqref{S3dagger:1} imply that
\begin{equation}\label{S3S4}
\begin{aligned}
\mathfrak{S}_{\delta, j}^{\Lambda*\mathbb{1}}
& \ll_{\varepsilon} 
\big((x^2 D^7)^{1/12}
+ x^{-1} D^2
+ DH^{-1}
+ (xD^{20}H)^{1/24}
+ (x^{-1}DH)^{1/2}
\big) x^{\varepsilon},
\end{aligned}
\end{equation}
for $j=3, 4$, $H\le D^{1/2}$ and $(L, M, N)$ verifying \eqref{Cond:LMN-Type2}.

\goodbreak

C. \textit{End of proof of \eqref{eq:prop_1}.}
\par
\vskip 0,5mm

Inserting \eqref{S1S2} and \eqref{S3S4} into \eqref{Decomposition} , we find that
\begin{align*}
\mathfrak{S}_{\delta}^{\Lambda*\mathbb{1}}(x, D)
& \ll \big((x^{171} D^{624})^{1/1026}
+ (x^2 D^7)^{1/12}
+ D^{9/10}
+ (x^{-1}D^3)^{1/2}
\\
& \hskip 3,5mm
+ DH^{-1}
+ (x^{19} D^{68} H^{19})^{1/114}
+ (x^5D^{17}H^5)^{1/30}
+ (xD^{20}H)^{1/24}
\big)x^{\varepsilon}
\end{align*}
for $H\le D^{20/57}$, 
where we have removed the term $(x^{-1}DH)^{1/2}\le (x^{-1}D^3)^{1/2}$.
Optimising $H$ over $[1, D^{20/57}]$, it follows that
\begin{align*}
\mathfrak{S}_{\delta}^{\Lambda*\mathbb{1}}(x, D)
& \ll \big((x^{171} D^{624})^{1/1026}
+ (x^2 D^7)^{1/12}
+ D^{9/10}
+ (x^{-1}D^3)^{1/2}
\\
& \hskip 3,5mm
+ (x^{19} D^{87})^{1/133}
+ (x^5 D^{22})^{1/35}
+ (xD^{21})^{1/25}
+ D^{9/19}
\\
& \hskip 3,5mm
+ (x^{19} D^{68})^{1/114}
+ (x^5D^{17})^{1/30}
+ (xD^{20})^{1/24}\big)x^{\varepsilon}
\end{align*}
which implies the required inequality, 
since the 5th term dominates the 7th one provided $D\le x^{171/308}$ 
and the first term dominates the others ones provided $D\le x^{57/103}$.
\end{proof}

The second proposition gives the required estimate for 
$\mathfrak{S}_{\delta}^{\tau_k}(x, D)$ defined as in \eqref{def:Stauk}.

\begin{proposition}\label{prop_2}
Let $k\ge 2$ be an integer.
For any $\varepsilon>0$, we have
\begin{equation}\label{prop_3:1}
\mathfrak{S}_{\delta}^{\tau_k}(x, D)
\ll_{k, \varepsilon} (x^kD^{4k-1})^{1/6k} x^{\varepsilon}
\end{equation}
uniformly for $x\ge 3$, $1\le D\le x^{\min\{4k/(5k+1), k/(2k-2)\}}$ and $0\le \delta\le \varepsilon^{-1}$,
where the implied constant depends on $k$ and $\varepsilon$ only.
\end{proposition}

\begin{proof}
Noticing that $\tau_k = \tau_{k-1}*\mathbb{1}$,
using Lemma \ref{lem:Vaaler} and splitting the interval of summation into the dyadic intervals, we can write
\begin{equation}\label{prop_3:2}
\begin{aligned}
\mathfrak{S}_{\delta}^{\tau_k}(x, D)
& = - \frac{1}{2\pi\text{i}} \sum_{H'} \sum_{M} \sum_{N}
\big(\mathfrak{S}_{\delta, \flat}^{\tau_k}(H', M, N) 
+ \overline{\mathfrak{S}_{\delta, \flat}^{\tau_k}(H', M, N)}\big)
\\
& \hskip 22mm
+ \sum_{M} \sum_{N} \mathfrak{S}_{\delta, \dagger}^{\tau_k}(M, N),
\end{aligned}
\end{equation}
where $H\le D$, $MN\asymp D$ (i.e. $D\ll MN\ll D$), $N\ge D^{1/k}$, $a_h:=\frac{H'}{h}\Phi\big(\frac{h}{H+1}\big)\ll 1$ and 
\begin{align*}
\mathfrak{S}_{\delta, \flat}^{\tau_k}(H', M, N)
& := \frac{1}{H'} \sum_{h\sim H'} a_h
\mathop{\sum_{m\sim M} \sum_{n\sim N}}_{D<mn\le 2D} \tau_{k-1}(m) \text{e}\Big(\frac{hx}{mn+\delta}\Big),
\\
\mathfrak{S}_{\delta, \dagger}^{\tau_k}(M, N)
& := \mathop{\sum_{m\sim M} \sum_{n\sim N}}_{D<mn\le 2D} 
\tau_{k-1}(m) R_H\Big(\frac{x}{mn+\delta}\Big).
\end{align*}

Firstly we bound $\mathfrak{S}_{\delta, \flat}^{\tau_k}(H', M, N)$.
Applying the exponent pair $(\frac{1}{2}, \frac{1}{2})$ to the sum over $n$,
it follows that 
\begin{equation}\label{prop_3:3}
\begin{aligned}
\mathfrak{S}_{\delta, \flat}^{\tau_k}(H', M, N)
& \ll \frac{x^{\varepsilon}}{H'} \sum_{h\sim H'} \sum_{m\sim M}
\Big\{\Big(\frac{xh}{mN^2}\Big)^{1/2} N^{1/2} + \Big(\frac{xh}{mN^2}\Big)^{-1}\Big\}
\\\noalign{\vskip 0,5mm}
& \ll \big((x H' M N^{-1})^{1/2} + x^{-1} D^2\big) x^{\varepsilon}.
\end{aligned}
\end{equation}
On the other hand, we remove the extra multiplicative condition $D<mn\le 2D$ at the cost of a factor $\log x$,
and then apply Lemma \ref{lem:LiuWuYang} with $\alpha=\beta=\gamma=1$,
$(X, H, M, N) = (xH'/MN, H', M, N)$ and $(\kappa, \lambda)=(\frac{1}{2}, \frac{1}{2})$ to get
$$
\mathfrak{S}_{\delta, \flat}^{\tau_k}(H', M, N)
\ll \big((xM^4N^3)^{1/6}
+ M N^{1/2}
+ H'^{-1/2} M^{1/2} N
+ (x^{-1}H'^{-1}M^3N^3)^{1/2}\big)X^{\varepsilon}
$$
provided $H'\le H\le D^{1/k}\;(\le N$).
Using $MN\asymp D$ and $N\ge D^{1/k}$, we can derive that
\begin{equation}\label{prop_3:4}
\begin{aligned}
\mathfrak{S}_{\delta, \flat}^{\tau_k}(H', M, N)
& \ll \big((x^kD^{4k-1})^{1/6k}
+ D^{1-1/(2k)}
+ H'^{-1/2} M^{1/2} N
+ (x^{-1}D^3)^{1/2}\big)x^{\varepsilon} 
\\
& \ll \big((x^kD^{4k-1})^{1/6k}
+ H'^{-1/2} M^{1/2} N\big)x^{\varepsilon} 
\end{aligned}
\end{equation}
for $H'\le H\le D^{1/k}$.
In the last inequality, we can remove $D^{1-1/(2k)}$ and $(x^{-1}D^3)^{1/2}$ 
because they can be absorbed by $(x^kD^{4k-1})^{1/6k}$ for $D\le x^{\min\{4k/(5k+1), k/(2k-2)\}}$.
From \eqref{prop_3:3} and \eqref{prop_3:4}, we deduce that
\begin{equation}\label{prop_3:5}
\mathfrak{S}_{\delta, \flat}^{\tau_k}(H', M, N)
\ll (x^kD^{4k-1})^{1/6k} x^{\varepsilon} 
\end{equation}
for $H'\le H\le D^{1/k}$ and $MN=D\le x^{\min\{4k/(5k+1), k/(2k-2)\}}$,
where we have used the following estimates
\begin{align*}
\min\{(x H' M N^{-1})^{\frac{1}{2}}, H'^{-\frac{1}{2}} M^{\frac{1}{2}} N\}
& \le ((x H' M N^{-1})^{\frac{1}{2}})^{\frac{1}{3}}
(H'^{-\frac{1}{2}} M^{\frac{1}{2}} N)^{\frac{2}{3}}
\\
& = (x H'^{-1} M^3 N^3)^{\frac{1}{6}}
\\
& \le (x D^3)^{1/6}
\end{align*}
and $\max\{(x D^3)^{1/6}, x^{-1} D^2\}\le (x^kD^{4k-1})^{1/6k}$ for $D\le x^{7k/(8k+1)}$.

On the other hand, \eqref{eq:Vaaler} of Lemma \ref{lem:Vaaler} allows us to derive that
\begin{align*}
\big|\mathfrak{S}_{\delta, \dagger}^{\tau_k}(M,N)\big|
& \le \sum_{m\sim M} \sum_{n\sim N} 
\Big|R_H\Big(\frac{x}{mn+\delta}\Big)\Big|
\\
& \le \frac{1}{2H+2} \sum_{0\le |h|\le H} \Big(1-\frac{|h|}{H+1}\Big)
\sum_{m\sim M} \sum_{n\sim N} {\rm e}\Big(\frac{xh}{d^2mn+\delta}\Big).
\end{align*}
When $h\not=0$, we can bound the triple sums as before.
Thus
\begin{equation}\label{prop_3:6}
\mathfrak{S}_{\delta, \dagger}^{\tau_k}(M, N)
\ll \big(DH^{-1} 
+ (x^kD^{4k-1})^{1/6k}\big) x^{\varepsilon}
\end{equation}
for $H\le D^{1/k}$ and $MN\asymp D\le x^{\min\{4k/(5k+1), k/(2k-2)\}}$.

Inserting \eqref{prop_3:5} and \eqref{prop_3:6} into \eqref{prop_3:2} and taking $H=D^{1/k}$, 
we can obtain the required inequality.
\end{proof}

The third proposition gives the required estimate for $\mathfrak{S}_{\delta}^{2^{\omega}}(x, D)$
 defined as in \eqref{def:S2omega}.

\begin{proposition}\label{prop_3}
For any $\varepsilon>0$, we have
$$
\mathfrak{S}_{\delta}^{2^{\omega}}(x, D)
\ll_{\varepsilon} (x^2 D^7)^{1/12} x^{\varepsilon}
$$
uniformly for $x\ge 3$ and $1\le D\le x^{8/11}$ and $0\le \delta\le \varepsilon^{-1}$,
where the implied constant depends on $\varepsilon$ only.
\end{proposition}

\begin{proof}
Firstly we write
\begin{align*}
\mathfrak{S}_{\delta}^{2^{\omega}}(x, D)
& = \sum_{d\le \sqrt{2D}} \mu(d) \sum_{D/d^2<\ell\le 2D/d^2} \tau(\ell) \psi\Big(\frac{x/d^2}{\ell+\delta/d^2}\Big)
\\
& = \sum_{d\le \sqrt{2D}} \mu(d) \mathfrak{S}_{\delta/d^2}^{\tau}(x/d^2, D/d^2).
\end{align*}
Secondly it is easy to see that $D\le x^{8/11}$ implies $D/d^2\le (x/d^2)^{8/11}$.
Thus we can apply Proposition \ref{prop_2} with $k=2$ to get the required result.
\end{proof}

\goodbreak

The last proposition gives the required estimate for $\mathfrak{S}_{\delta}^{\mu_2}(x, D)$
defined as in \eqref{def:Smu2}.

\begin{proposition}\label{prop_4}
For any $\varepsilon>0$, we have
$$
\mathfrak{S}_{\delta}^{\mu_2}(x, D)
\ll_{\varepsilon} \big((xD^3)^{1/7} + (xD^2)^{1/6}\big)x^{\varepsilon}
$$
uniformly for $x\ge 3$, $1\le D\le x^{8/11}$ and $0\le \delta\le \varepsilon^{-1}$,
where the implied constant depends on $\varepsilon$ only.
\end{proposition}

\begin{proof}
Using Lemma \ref{lem:Vaaler} and splitting the interval of summation $[1, 2D]$ 
into the dyadic intervals $(L, 2L]$, we can write
\begin{equation}\label{proof:prop_2:1}
\mathfrak{S}_{\delta}^{\mu_2}
= - \frac{1}{2\pi\text{i}} \sum_{L}
\big(\mathfrak{S}_{\delta, \flat}^{\mu_2}(L) 
+ \overline{\mathfrak{S}_{\delta, \flat}^{\mu_2}(L)}\big)
+ \sum_{L} \mathfrak{S}_{\delta, \dagger}^{\mu_2}(L),
\end{equation}
where $H\le D$ and 
\begin{align*}
\mathfrak{S}_{\delta, \flat}^{\mu_2}(L)
& := \sum_{h\le H} \frac{1}{h}\Phi\big(\frac{h}{H+1}\big)
\sum_{\ell\sim L} \sum_{(D/\ell)^{1/2}<d\le (2D/\ell)^{1/2}} 
\mu(d) \text{e}\Big(\frac{hx}{d^2\ell+\delta}\Big),
\\
\mathfrak{S}_{\delta, \dagger}^{\mu_2}(L)
& := \sum_{\ell\sim L} \sum_{(D/\ell)^{1/2}<d\le (2D/\ell)^{1/2}}  
\mu(d) R_H\Big(\frac{x}{d^2\ell+\delta}\Big).
\end{align*}
Inverting the order of summations and applying the exponent pair $(\frac{1}{6}, \frac{4}{6})$ to the sum over $\ell$,
it follows that
\begin{equation}\label{proof:prop_2:2}
\begin{aligned}
\mathfrak{S}_{\delta, \flat}^{\mu_2}(L)
& \ll \sum_{h\sim H} \frac{1}{h} \sum_{(D/2L)^{1/2}<d\le (2D/L)^{1/2}}
\bigg(\Big(\frac{xh}{d^2L^2}\Big)^{1/6} L^{4/6} + \Big(\frac{xh}{d^2L^2}\Big)^{-1}\bigg)
\\\noalign{\vskip 1mm}
& \ll \big((x D^2 H)^{1/6} + x^{-1} D^2\big) x^{\varepsilon}.
\end{aligned}
\end{equation}
On the other hand, \eqref{eq:Vaaler} of Lemma \ref{lem:Vaaler} allows us to derive that
\begin{align*}
\big|\mathfrak{S}_{\delta, \dagger}^{\mu_2}(L)\big|
& \le \sum_{(D/2L)^{1/2}<d\le (2D/L)^{1/2}} \sum_{\ell\sim L} 
\Big|R_H\Big(\frac{x}{d^2\ell+\delta}\Big)\Big|
\\
& \le \frac{1}{2H+2} \sum_{0\le |h|\le H} \Big(1-\frac{|h|}{H+1}\Big) 
\sum_{(D/2L)^{1/2}<d\le (2D/L)^{1/2}} \sum_{\ell\sim L} {\rm e}\Big(\frac{xh}{d^2\ell+\delta}\Big).
\end{align*}
When $h\not=0$, as before we apply the exponent pair $(\frac{1}{6}, \frac{4}{6})$ to the sum over $\ell$ and obtain
\begin{equation}\label{proof:prop_2:3}
\mathfrak{S}_{\delta, \dagger}^{\mu_2}(L)
\ll \big(DH^{-1} + (x D^2 H)^{1/6} + x^{-1} D^2\big) x^{\varepsilon}.
\end{equation}
Inserting \eqref{proof:prop_2:2} and \eqref{proof:prop_2:3} into \eqref{proof:prop_2:1}, it follows that
$$
\mathfrak{S}_{\delta}^{\mu_2}
\ll \big(DH^{-1} + (x D^2 H)^{1/6} + x^{-1} D^2\big) x^{\varepsilon}
$$
for $H\le D$. Optimising $H$ on $[1, D]$, we find that
$$
\mathfrak{S}_{\delta}^{\mu_2}
\ll \big((x D^3)^{1/7} + (x D^2)^{1/6} + x^{-1} D^2\big) x^{\varepsilon}
$$
This implies the required result, since $x^{-1} D^2\le (x D^3)^{1/7}$ for $D\le x^{8/11}$.
\end{proof}

\vskip 5mm

\section{Proof of Theorem \ref{thm}}

Let $f=\omega$ or $\tau_k$ or $2^{\omega}$ or $\mu_2$
and let $N_f\in [1, x^{1/2})$ be a parameter which can be chosen later.
First we write
\begin{equation}\label{4.1}
S_f(x) = \sum_{n\le x} f\Big(\Big[\frac{x}{n}\Big]\Big) = S_f^{\dagger}(x)+S_f^{\sharp}(x)
\end{equation}
with
$$
S_f^{\dagger}(x):=\sum_{n\le N_f} f\Big(\Big[\frac{x}{n}\Big]\Big),      
\qquad 
S_f^{\sharp}(x):=\sum_{N_f<n\le x} f\Big(\Big[\frac{x}{n}\Big]\Big).
$$
In our case, we have $f(n)\ll_{\varepsilon} n^{\varepsilon}$ for any $\varepsilon>0$ and all $n\ge 1$.
Thus
\begin{equation}\label{4.2}
S_f^{\dagger}(x)
\ll_{\varepsilon} N_f x^{\varepsilon}.
\end{equation}
In order to bound $S_f^{\sharp}(x)$, we put $d=[x/n]$.
Noticing that
$$
x/n-1<d\le x/n
\;\Leftrightarrow\;
x/(d+1)<n\le x/d,
$$
we can derive that
\begin{equation}\label{4.3}
\begin{aligned}
S_f^{\sharp}(x)
& = \sum_{d\le x/N_f} f(d) \sum_{x/(d+1)<n\le x/d} 1
\\
& = \sum_{d\le x/N_f} f(d) \Big(\frac{x}{d}-\psi\Big(\frac{x}{d}\Big)-\frac{x}{d+1}+\psi\Big(\frac{x}{d+1} \Big) \Big)
\\
& = x \sum_{d\ge 1} \frac{f(d)}{d(d+1)} + \mathcal{R}_1^{f}(x, N_f) - \mathcal{R}_0^{f}(x, N_f) + O(N_f),
\end{aligned}
\end{equation}
where we have used the following bounds
$$
x \sum_{d>x/N_f} \frac{f(d)}{d(d+1)}\ll_{\varepsilon} N_f x^{\varepsilon},
\qquad
\sum_{d\le N_f} f(d)\Big(\psi\Big(\frac{x}{d+1}\Big) - \psi\Big(\frac{x}{d}\Big)\Big)
\ll_{\varepsilon} N_f x^{\varepsilon}
$$
and
$$
\mathcal{R}_{\delta}^{f}(x, N_f)
= \sum_{N_f<d\le x/N_f} f(d) \psi\Big(\frac{x}{d+\delta}\Big).
$$
Combining \eqref{4.1}, \eqref{4.2} and \eqref{4.3}, it follows that
$$
S_f(x) = x \sum_{d\ge 1} \frac{f(d)}{d(d+1)} 
+ O_{\varepsilon}\big(|\mathcal{R}_1^{f}(x, N_f)| + |\mathcal{R}_0^{f}(x, N_f)| + N_fx^{\varepsilon}\big).
$$
Thus in order to prove Theorem \ref{thm}, it suffices to show that
\begin{equation}\label{UB:Rdeltaf}
\mathcal{R}_{\delta}^{f}(x, N_f)\ll N_f x^{\varepsilon}
\qquad
(x\ge 1)
\end{equation}
for 
$$\displaylines{
\vbox{\tabskip = 0mm\offinterlineskip
\halign{
\vrule # & &\hfil$ $ $#$  \hfil & \vrule #\cr
\noalign{\hrule}
height 1,5mm
&&&&&&&&&&
\cr
& f & 
& \omega & 
& \tau_k & 
& 2^{\omega} & 
& \mu_2 & 
\cr
height 1,5mm
&&&&&&&&&&
\cr
\noalign{\hrule}
height 1,5mm
&&&&&&&&&&
\cr
& \;\; N_f \;\; &
& \;\; x^{53/110} \;\; & 
& \;\; x^{(5k-1)/(10k-1)} \;\; &
& \;\; x^{9/19} \;\; & 
& \;\; x^{2/5} \;\; & 
\cr
height 1,5mm
&&&&&&&&&&
\cr
\noalign {\hrule}}}
\cr}$$

\subsection{Proof of \eqref{UB:Rdeltaf} with $f=\omega$}\

\vskip 0,5mm

We have $N_{\omega}=x^{53/110}$.
In order to apply \eqref{eq:prop_1} of Proposition \ref{prop_1}, we must switch from
$$
\mathcal{R}_{\delta}^{\omega}(x, N_{\omega})
= \sum_{N_{\omega}<p\ell\le x/N_{\omega}} \psi\Big(\frac{x}{p\ell+\delta}\Big) 
= \sum_{N_{\omega}<p\ell\le x/N_{\omega}} \frac{\Lambda(p)}{\log p} \psi\Big(\frac{x}{p\ell+\delta}\Big) 
$$
to 
$$
\widetilde{\mathcal{R}}_{\delta}^{\omega}(x, N_{\omega}, t)
:= \sum_{d\le t} \Lambda(d) \sum_{N_{\omega}/d<\ell\le (x/N_{\omega})/d} \psi\Big(\frac{x}{d\ell+\delta}\Big).
$$
For this, we need to estimate the contribution of prime powers:
\begin{align*}
\mathscr{R}_1
& := \sum_{\substack{N_{\omega}<p^{\nu}\ell\le x/N_{\omega}\\ \ell\le (x/N_{\omega})^{53/114}, \, \nu\ge 2}} \frac{\Lambda(p^{\nu})}{\log p^{\nu}} \psi\Big(\frac{x}{p^{\nu}\ell+\delta}\Big)
= \sum_{\substack{N_{\omega}<p^{\nu}\ell\le x/N_{\omega}\\ \ell\le (x/N_{\omega})^{53/114}, \, \nu\ge 2}} 
\frac{1}{\nu} \psi\Big(\frac{x}{p^{\nu}\ell+\delta}\Big),
\\
\mathscr{R}_2
& := \sum_{\substack{N_{\omega}<p^{\nu}\ell\le x/N_{\omega}\\ \ell>(x/N_{\omega})^{53/114}, \, \nu\ge 2}} \frac{\Lambda(p^{\nu})}{\log p^{\nu}} \psi\Big(\frac{x}{p^{\nu}\ell+\delta}\Big)
= \sum_{\substack{N_{\omega}<p^{\nu}\ell\le x/N_{\omega}\\ \ell>(x/N_{\omega})^{53/114}, \, \nu\ge 2}} 
\frac{1}{\nu} \psi\Big(\frac{x}{p^{\nu}\ell+\delta}\Big).
\end{align*}
Firstly we have trivially
$$
\mathscr{R}_1
\ll \sum_{\ell\le (x/N_{\omega})^{53/114}} (x/N_{\omega}\ell)^{1/2}
\ll (x/N_{\omega})^{167/228}
\ll x^{167/440}.
$$
On the other hand, applying \cite[Lemma 2.4]{Zhai2020} with $(\kappa, \lambda)=(\frac{1}{2}, \frac{1}{2})$, 
we can derive that
\begin{align*}
\mathscr{R}_2
& = \sum_{\substack{p^{\nu}\le (x/N_{\omega})^{61/114}\\ \nu\ge 2}} \frac{1}{\nu} 
\sum_{(x/N_{\omega})^{53/114}<\ell\le x/N_{\omega}p^{\nu}} \psi\Big(\frac{x}{p^{\nu}\ell+\delta}\Big)
\\
& \ll \sum_{\substack{p^{\nu}\le (x/N_{\omega})^{61/114}\\ \nu\ge 2}} 
\max_{\substack{(x/N_{\omega})^{53/114}<N_1\le x/N_{\omega}p^{\nu}\\ N_1<N_2\le 2N_1}} 
N_1\bigg|\sum_{N_1<\ell\le N_2} \frac{1}{\ell} \psi\Big(\frac{x}{p^{\nu}\ell+\delta}\Big)\bigg|
\\
& \ll \sum_{\substack{p^{\nu}\le (x/N_{\omega})^{61/114}\\ \nu\ge 2}} 
\max_{\substack{(x/N_{\omega})^{53/114}<N_1\le x/N_{\omega}p^{\nu}\\ N_1<N_2\le 2N_1}} 
\bigg(\frac{N_1^2}{(x/p^{\nu})} + (x/p^{\nu})^{1/3}\bigg)
\\
& \ll \sum_{\substack{p^{\nu}\le (x/N_{\omega})^{61/114}\\ \nu\ge 2}} 
\bigg(\frac{x}{N_{\omega}^2p^{\nu}} + (x/p^{\nu})^{1/3}\bigg)
\\
& \ll x/N_{\omega}^2 + x^{1/3} (x/N_{\omega})^{(61/114)/6}
\\\noalign{\vskip 2mm}
& \ll x^{167/440}
\end{align*}
Thus we can write
\begin{equation}\label{R-tildeR}
\begin{aligned}
\mathcal{R}_{\delta}^{\omega}(x, N_{\omega})
& = \sum_{N_{\omega}<d\ell\le x/N_{\omega}} \frac{\Lambda(d)}{\log d} \psi\Big(\frac{x}{d\ell+\delta}\Big) 
+ O\big(x^{167/440}\big)
\\
& = \int_{2-}^{x/N_{\omega}} \frac{1}{\log t} \d\widetilde{\mathcal{R}}_{\delta}^{\omega}(x, N_{\omega}, t) 
+ O\big(x^{167/440}\big)
\\\noalign{\vskip 1,5mm}
& \ll \max_{2\le t\le x/N_{\omega}} \big|\widetilde{\mathcal{R}}_{\delta}^{\omega}(x, N_{\omega}, t)\big| + x^{167/440}.
\end{aligned}
\end{equation}
Writing $D_j:=x/(2^jN_{\omega})$, we have $N_{\omega}\le D_j\le x/N_{\omega}\le x^{57/110}$ 
for $0\le j\le \frac{\log(x/N_{\omega}^2)}{\log 2}$.
Thus we can apply \eqref{eq:prop_1} of Proposition \ref{prop_1} to get
\begin{align*}
|\widetilde{\mathcal{R}}_{\delta}^{\omega}(x, N_{\omega}, t)|
& \le \sum_{0\le j\le \log(x/N_{\omega}^2)/\log 2} |\mathfrak{S}_{\delta}^{\Lambda*\mathbb{1}}(x, D_j)|
\\
& \ll \sum_{0\le j\le \log(x/N_{\omega}^2)/\log 2} \big((x^{19} D_j^{87})^{1/133} + (x^{171} D_j^{624})^{1/1026}\big) x^{\varepsilon}
\\
& \ll \big((x^{106} N_{\omega}^{-87})^{1/133} + (x^{795} N_{\omega}^{-624})^{1/1026}\big) x^{\varepsilon}
\end{align*}
for all $2\le t\le x/N_{\omega}$.
Inserting this into \eqref{R-tildeR}, we get the required result.

\subsection{Proof of \eqref{UB:Rdeltaf} with $f=\tau_k$}\

\vskip 0,5mm

We have $N_{\tau_k}=x^{(5k-1)/(10k-1)}$.
Let $D_j:=x/(2^jN_{\tau_k})$, then we easily see that 
$N_{\tau_k}\le D_j\le x/N_{\tau_k}\le x^{\min\{4k/(5k+1), k/(2k-2)\}}$ 
for $0\le j\le \frac{\log(x/N_{\tau_k}^2)}{\log 2}$.
Thus Proposition \ref{prop_2} gives us
\begin{align*}
|\mathcal{R}_{\delta}^{\tau_k}(x, N_{\tau_k})|
& \le \sum_{0\le j\le \log(x/N_{\tau_k}^2)/\log 2} |\mathfrak{S}_{\delta}^{\tau_k}(x, D_j)|
\\
& \ll \sum_{0\le j\le \log(x/N_{\tau_k}^2)/\log 2} (x^kD_j^{4k-1})^{1/6k} x^{\varepsilon}
\\
& \ll (x^{5k-1}N_{\tau_k}^{-(4k-1)})^{1/6k} x^{\varepsilon},
\end{align*}
which implies the required result.

\subsection{Proof of \eqref{UB:Rdeltaf} with $f=2^{\omega}$}\

\vskip 0,5mm

We have $N_{2^{\omega}}=x^{9/19}$.
Let $D_j:=x/(2^jN_{2^{\omega}})$, then $N_{2^{\omega}}\le D_j\le x/N_{2^{\omega}}\le x^{10/19}$ 
for $0\le j\le \frac{\log(x/N_{2^{\omega}}^2)}{\log 2}$.
Noticing that $2^{\omega} = \tilde{\mu}*\tau$, we can apply Proposition \ref{prop_3} to get
\begin{align*}
|\mathcal{R}_{\delta}^{2^{\omega}}(x, N_{2^{\omega}})|
& \le \sum_{0\le j\le \log(x/N_{2^{\omega}}^2)/\log 2} |\mathfrak{S}_{\delta}^{2^{\omega}}(x, D_j)|
\\
& \ll \sum_{0\le j\le \log(x/N_{2^{\omega}}^2)/\log 2} (x^2D_j^7)^{1/12} x^{\varepsilon} 
\\
& \ll (x^9N_{2^{\omega}}^{-7})^{1/12} x^{\varepsilon},
\end{align*}
which implies the required result.

\subsection{Proof of \eqref{UB:Rdeltaf} with $f=\mu_2$}\

\vskip 0,5mm

We have $N_{\mu_2}=x^{2/5}$.
Let $D_j:=x/(2^jN_{\mu_2})$, then $N_{\mu_2}\le D_j\le x/N_{\mu_2}\le x^{3/5}$ 
for $0\le j\le \frac{\log(x/N^2_{\mu_2})}{\log 2}$.
Noticing that $\mu_2 = \tilde{\mu}*\mathbb{1}$, we can apply Proposition \ref{prop_4} to get
\begin{align*}
|\mathcal{R}_{\delta}^{\mu_2}(x, N_{\mu_2})|
& \le \sum_{0\le j\le \log(x/N_{\mu_2}^2)/\log 2} |\mathfrak{S}_{\delta}^{\mu_2}(x, D_j)|
\\
& \ll \sum_{0\le j\le \log(x/N_{\mu_2}^2)/\log 2} \big((xD_j^3)^{1/7} + (xD_j^2)^{1/6}\big)x^{\varepsilon}
\\
& \ll \big((x^4 N_{\mu_2}^{-3})^{1/7} + (x^3 N_{\mu_2}^{-2})^{1/6}\big) x^{\varepsilon},
\end{align*}
which implies the required result.

\vskip 3mm

\noindent{\bf Acknowledgement}.
This work is supported in part by the National Natural Science Foundation of China 
(Grant Nos. 12071238, 11771252, 11971370 and 12071375)
and by the NSF of Chongqing (Grant  No. cstc2019jcyj-msxm1651).

\end{document}